\documentclass[12pt]{article}
\usepackage{amssymb, amsmath, amsfonts}
\def\nodo#1{}
\def\genfd{{\bf k}}
\def\widem{\widehat{\mathcal M}}

\font\tenDDl=msbm10  
\font\sevenDDl=msbm7 
\font\fiveDDl=msbm5 
        \newfam\DDlfam \def\DDl{\fam\DDlfam\tenDDl} 
                \textfont\DDlfam=\tenDDl
\scriptfont\DDlfam=\sevenDDl
        \scriptscriptfont\DDlfam=\fiveDDl
\def\nxpoint{\refstepcounter{section}%
  \makepoint{\thesection}}
\def\refpt#1{{\rm\textbf{\ref{#1}}}}
\def\makepoint#1{\medbreak\noindent{\bf #1. }}

\begin{document}
\begin{center}
{\bf \Large A simple algorithm for extending the identities
for quantum minors to the multiparametric case}
\vskip .1in
\centerline{Zoran \v{S}koda, {\tt zskoda@irb.hr}}
\vskip .1in
{\bf Abstract.} {\it For any homogeneous identity between $q$-minors,
we provide an identity between $P,Q$-minors.
}
\end{center}

\nxpoint (\cite{Artin:multi, demidov:rev, Manin:multi}) 
Suppose we are given $n\times n$ matrices $P = (p_{ij})$ and 
$Q = (q_{ij})$ with invertible entries in the ground field $\genfd$, 
for which there exist $q$ such that  
\begin{equation}\label{eq:pqmu}
  p_{ij} q_{ij} = q^2, q_{ij} = q^{-1}_{ji},\,\,\,\,\,i < j,\mbox{ and }\,\,
q_{ii} = p_{ii},\,\,\, \mbox{for all }i. \end{equation}

Define an associative algebra
${\cal M}(P,Q;\genfd,n):= \genfd \langle T^i_j , i,j = 1,\ldots, n\rangle/I$,
where $I$ is the ideal spanned by the relations
\[\begin{array}{ll}
T^k_i T^k_j = q_{ij} T^k_j T^k_i, & i < j \\
T^k_i T^l_i = p_{kl} T^l_i T^k_i, & k < l \\
q_{ij} T^k_j T^l_i = p_{kl} T^l_i  T^k_j, & i<j,k<l \\
T^k_i T^l_j - q_{ij} q^{-1}_{kl} T^l_j T^k_i
= (q_{ij}-p_{ij}^{-1}) T^k_j T^l_i,& i<j,k<l
\end{array}\]
${\mathcal M} := {\mathcal M}(P,Q;\genfd,n)$ is a bialgebra 
with respect to the "matrix" comultiplication
which is the unique algebra homomorphism $\Delta : {\mathcal M} 
\to {\mathcal M} \otimes {\mathcal M}$
extending the formulas $\Delta T^i_j = \sum T^i_k \otimes T^k_j$ with 
counit $\epsilon T^i_j = \delta^i_j$ (Kronecker delta). 
The bialgebra $\mathcal M$ is called 
the {\bf multiparametric quantum linear group}. 

Our conventions differ a bit from the cited references: 
we treat $p$-s and $q$-s symmetrically in the sense
that if we interchange rows and columns of matrix $T$ 
and if we simultaneously interchange $P$ and $Q$, 
we obtain an isomorphic algebra.
If $P = Q$ and
$q_{ij} = q$ for $i < j$ and $q_{ij} = q^{-1}$ for $i > j$, then
${\mathcal M} = {\mathcal M}_q(\genfd)$ 
(1-parametric quantized matrix bialgebra).
In this paper, we will denote by $t^i_j$ the generators for 1-parametric case.

\nxpoint ({\it Labels.}) 
It is convenient to consider that the row and column labels
belong to some totally ordered set of labels, 
not necessarily the set $\{ 1, \ldots, n\}$. 
The main reason is that one often needs to treat some
subsets of the set of labels and 
the corresponding submatrices of the matrix $T = (T^i_j)$.

\nxpoint The quantum $Q$-space $\mathcal{O}(\genfd_Q^n)$,
quantum $P$-space $\mathcal{O}(\genfd_P^n)$, 
$q$-{ normalized} $Q$-{ space}
$S_r(q,Q)$, { dual} $q$-{ normalized} $Q$-{ space}
$S_l(q,Q)$, right $P$-exterior algebra $\Lambda_P$, 
and left $Q$-exterior algebra $\Lambda_Q$, 
are the algebras {\it defined} by generators and relations
as follows:
\begin{equation}\label{eq:ordalg}\begin{array}{l}
\mathcal{O}(\genfd_Q^n) := \genfd \langle x^i, i = 1,\ldots, n \rangle /
\langle x^i x^j - q_{ij} x^j x^i, i < j \rangle \\
\mathcal{O}(\genfd_P^n) := \genfd \langle y_i, i = 1,\ldots, n \rangle /
\langle y_i y_j - p_{ij} y_j y_i, i < j \rangle \rangle \\
S_r(q,Q) := \genfd \langle r_1,\ldots, r_n \rangle/
\langle r_i r_j - qq_{ij}^{-1} r_j r_i, \,\,\,i<j\rangle \\
S_l(q,Q) := \genfd \langle l^1,\ldots, l^n \rangle/
\langle l^i l^j - q_{ij}q^{-1} l^j l^i,\,\,\,i<j\rangle \\
\Lambda_P := \genfd \langle e^1, \ldots, e^n \rangle /
\langle e^i e^j + p_{ij}^{-1} e^j e^i , \,i<j,\, (e^i)^2 \rangle \\
\Lambda_Q := \genfd \langle e_1, \ldots, e_n \rangle /
\langle f_i f_j + q_{ij}^{-1} f_j f_i , \,i<j,\, (f_i)^2 \rangle \\
\end{array}\end{equation}
We note after {\sc Manin} (e.g. \cite{ManinQGNG,Manin:multi}) 
the simple fact that 
$\mathcal{O}(\genfd_Q^n), \Lambda_Q$ are right and
$\mathcal{O}(\genfd_P^n), \Lambda_P$ are left comodule algebras
over $\mathcal M$ via coactions $\rho, \rho_\Lambda, \rho', \rho'_\Lambda$
which are the unique coactions which are algebra maps and extend formulas
$\rho(x^i) = \sum_j x^j \otimes T^i_j$, $\rho_\Lambda(e^i) = 
\sum_j e^j\otimes T^i_j$, $\rho'(y_j)=\sum_i T^i_j\otimes y_i$ and
$\rho'_\Lambda(f_j) =\sum_i T^i_j\otimes f_i$.

\nxpoint Tuples of labels will be called {\it multilabels}.  
Let $L = (l_1,\ldots, l_r)$, $K = (k_1,\ldots,k_s)$. 
The concatenation will be denoted by juxtaposition:
$LK = (l_1,\ldots, l_r, k_1,\ldots, k_s)$. Usually the multilabels will be 
(ascendingly) ordered to start with and $\hat{L}$ 
denotes the ordered complement of a submultilabel
$L$ (usually in $\{1,\ldots, n\}$). 
By placing the multilabel within the colons, 
we will denote its ascendingly ordered version.
For example, if $K$ and $L$ are ordered, then $KL$ is not necessarily ordered,
because some labels in $L$ may be smaller than some labels in $K$. 
However $:KL:$ is the multilabel obtained 
from $KL$ by permuting the labels until they are ascendingly
ordered. We identify the notation 
for a single label $j$ and the multilabel $(j)$,
and in this vein $\hat{j} = \widehat{(j)}$ is the same as multilabel 
$(1,2,\ldots, j-1, j+1,\ldots, n)$.

We use obvious exponent notation: $r^J := r^{j_1}\cdots r^{j_n}$ and alike.

\nxpoint Algebras~(\ref{eq:ordalg}) satisfy the obvious normal basis 
(``PBW type'') theorems: fix an order on 
generators, then the monomials ordered compatibly with this order
form a vector space basis (``PBW basis''), for example 
$x_1^{a_1} x_2^{a_2}\cdots x_n^{a_n}$ 
in the case of $\mathcal{O}(\genfd_Q^n)$; however, no
higher exponents than $1$ appear in the bases 
for $\Lambda_P$ and $\Lambda_Q$ because $e_i^2 = 0$.

Passing from arbitrary words in generators to the basis elements
clearly reduces to reordering the generators, and accumulating
the proportionality constants. Let us introduce the notation for 
those rearrangements factors which will be needed below. The following are
the defining properties of coefficient functions 
$\epsilon_P, \epsilon_Q,\zeta_r,\zeta_l$ from sets of
multilabels to $\genfd$:
$$\begin{array}{ll}
e_J := \epsilon_P(J) e_{:J:},&r_J := \zeta_r(J) r_{:J:} ,\\
f_J := \epsilon_Q(J) f_{:J:},&l_J := \zeta_l(J) l_{:J:} ,\\
\end{array}$$
where $\epsilon_P(J),\epsilon_Q(J)$ are defined only when $J$ has no
repeted labels inside, but
$\zeta_r$ and $\zeta_l$ are defined even for multilabels with repetition.
The following is obvious:
$$
\zeta_r(k_1, \ldots, k_s) = \prod_{i<j, k_i > k_j} q^{-1}q_{k_i k_j} 
$$
$$
\zeta_l(k_1, \ldots, k_s) = \prod_{i<j, k_i > k_j} qq^{-1}_{k_i k_j}
= \zeta_r^{-1}(k_1,\ldots,k_s)
$$
Clearly, if $J$ is ascendingly ordered multilabel with out repetitions, then
$$\begin{array}{lcccl}
\epsilon_P(\sigma J)&=& (-q)^{-l(\sigma)} \zeta_r(\sigma J)
&= &\prod_{k<l, \sigma (k) > \sigma(l)} (-p_{j_{\sigma(k)} j_{\sigma(l)}}),\\
\epsilon_Q(\sigma J)&=& (-q)^{l(\sigma)} \zeta_l(\sigma J)
&= &\prod_{k<l, \sigma (k) > \sigma(l)}(-q_{j_{\sigma(k)} j_{\sigma(l)}}),
\end{array}$$
(recall that $(q_{ij})^{-1} = q_{ji}$).

\nxpoint ({\it Gradings.})
Let ${\mathcal I}$ be the set of labels of generators
of $S_r = S_r(q,Q)$ (they label rows!).
Let ${\mathcal J}$ be the set of labels of generators
of $S_l = S_l(q,Q)$ (they label columns!).
Both sets are bijective to $\{1,\ldots,n\}$.
Thus the free Abelian group ${\DDl Z}[{\mathcal I}]$
is isomorphic to ${\DDl Z}^n$ (and could be naturally identified
with the weight lattice for $SL_n$).
Now we assign
${\DDl Z}[{\mathcal I}]{\rm -}{\DDl Z}[{\mathcal J}]$
bigrading to algebras $S_l(q,Q)$, $S_r(q,Q)$ and ${\cal M}_q(\genfd)$.
If $i_1,\ldots,i_n,j_1,\ldots,j_n$ are the elements of $\mathcal I$ and
$\mathcal J$, then a bidegree is a formal sum of the form 
$a_1 i_1 + \ldots + a_n i_n + b_1 j_1 +
\ldots b_n j_n$, e.g. $-i_3 + 2i_4 + j_3$, 
and we may separate the $\mathcal I$ and $\mathcal J$ grading
with comma for clarity, 
e.g. $(-i_3+2i_4,\, j_3) \equiv -i_3 + 2i_4 + j_3$.
We assign the bidegree $(-i,0)$
to the generator $r_i$ of $S_l$ (notice the negative sign!)
${\DDl Z}\langle {\cal I}^*\rangle$-degree to zero
and similarly the dual prescription $(0,-j)$ to $l^j$ of $S_r$.
We also assign the bidegree $(+i,+j)$ to each generator
$t^i_j$ of the 1-parametric algebra ${\cal M}_q(\genfd)$. 
The defining ideals are bihomogeneous hence we extend 
this prescription multiplicatively to a bigrading 
on the algebras $S_l(q,Q)$, $S_r(q,Q)$ and ${\cal M}_q(\genfd)$.

\nxpoint {\bf Notation.} Consider the tensor product of bigraded algebras
\begin{equation}\label{eq:big-algebra-MSS}
\widem := S_l(q,Q) \otimes {\cal M}_q(\genfd) \otimes S_r(q,Q).
\end{equation}
\nxpoint \label{lem:tensor1}{\bf Lemma.} {\it Any (bi)homogenous element in
a tensor product of (bi)graded algebras is a sum of
tensor products of (bi)homogenous elements in tensor factors. If one chooses
a set of homogeneous generators in each tensor factor 
than the summands can be chosen as tensor products of monomials 
in those generators. }

\nxpoint \label{lem:tensor2} {\bf Lemma.} {\it $l^j \otimes t^i_j \otimes r_i$ 
generate the subalgebra of all elements of bidegree $(0,0)$ in $\widem$.
}

{\it Proof.} If $J = (j_1,\ldots,j_s)$ is some ordered
$s$-tuple of labels (repetitions of labels possible)
denote $l^J = l^{j_1}\cdots l^{j_s}$ and we adopt obvious 
extension of this multilabel notation for $r$-s and $t$-s.
It is clear that any tensor product of monomials which is
of bidegree $(0,0)$ is of the form 
$l^{\tau I}\otimes t^{J}_{I}\otimes r_{\sigma J}$
where $\sigma$ and $\tau$ are permutations on $|I|=|J|$ letters.
Then by lemma~\refpt{lem:tensor1} 
it is enough to show that
any such tensor product $l^{\tau I}\otimes t^{J}_{I}\otimes r_{\sigma J}$ 
may be written as sum of products of the form $l^j \otimes t^i_j \otimes r_i$.
But $l^J$ and $l^{\tau J}$ are proportional in $S_l$, and similarly
$r_I$ and $r_{\sigma I}$ are proportional in $S_r$. Hence, up to accounting for
a scalar factor, we may assume that $\sigma$ and $\tau$ are trivial. But 
then the expression is manifestly the product of elements of the required form.

\nxpoint {\bf Theorem.} {\it Suppose ~(\ref{eq:pqmu}) holds.
Let $t^i_j$ and $T^i_j$ denote the generators of $q$-deformed
and $P,Q$-deformed quantum matrix algebras respectively. Then

(i) the rule 
\begin{equation}\label{eq:iota_SLq_to_multipar}
 \iota_{q,Q} : T^i_j \mapsto l^j \otimes t^i_j \otimes r_i
\end{equation}
extends to a unique algebra homomorphism 
$\iota_{q,Q} : {\cal M}(P,Q;\genfd)\to\widem$.  

(ii) This homomorphism is injective and
its image is the subalgebra of all elements 
of $(0,0)$-bidegree in $\widem$.

(iii) Similarly, rescaling $e^j$ by $l^j$ produces 
the relations in $\Lambda_P$ from
the relations in $\Lambda_q$.

}

{\it Proof}. (i) One needs to show that $\iota_{q,Q}$ sends
the ideal of relations (in free algebra on $T$-s) to zero.
For example, omitting the tensor product notation, we calculate,
for $i<j$ and $k<l$,\[\begin{array}{lcl}
\iota(q_{kl}T^i_l T^j_k) &=& q_{kl} l^l  t^i_l r_i \cdot l^k t^j_k r_j\\
&=& q_{kl} (l^l l^k)  (t^i_l t^j_k)  (r_i r_j) \\
&=& q_{kl} (q q_{kl}^{-1} l^k l^l) (t^j_k t^i_l) (q q_{ij}^{-1} r_j r_i)
\\
&=& q^2 q_{ij}^{-1} l^k l^l t^j_k t^i_l r_j r_i \\
&=& \iota(p_{ij}T^j_k T^i_l). \end{array}\]
 The other cases are left to the reader.

(ii) For injectivity one can use e.g. the normal basis theorem for the
quantum matrix algebras: monomials of the form $(T^1_1)^{\alpha_{11}}
(T^1_1)^{\alpha_{12}}\cdots(T^1_1)^{\alpha_{nn}}$ make a basis of 
${\mathcal M}(P,Q;\genfd)$. It is clear that the images are linearly 
independent because the middle tensor factors of the images are such (by
the normal basis theorem for 1-parametric case) and
the other two tensor factors are nonzero. The description of the image
of $\iota_{q,Q}$ follows from \refpt{lem:tensor2}.

(iii) Easy.

\nxpoint {\bf Remarks.}  This isomorphism will be very useful for our purpose. 
Essentially this proposition is a mechanism essentially equivalent  
to the cocycle-twisting of~\cite{Artin:multi}.
Namely, in both approaches, the difference 
between the algebra relations for  $t^i_j$-s and
for $T^i_j$-s is reflected in rescaling factors for each monomial, 
which may be expressed in terms of a bicharacter and 
depends only on the bidegree of the monomial. 

However, there is an important difference in using our isomorphism
$\iota$ from the usage of twisting in~\cite{Artin:multi}. 
Namely, it is shown in~\cite{Artin:multi}
that the correspondence $t^i_j \mapsto T^i_j$ which they use, 
extends multiplicatively on monomials 
to an isomorphisms of vector spaces, and even of coalgebras; whereas 
it does not respect the algebra structure. On the other hand, our map $\iota$
is a monomorphism of algebras, as stated above, but it does not respect the
coalgebra structure! 

\nxpoint ({\it Quantum minors.}) If $J$ and $K$ are ascendingly ordered
row and column multilabels of the same cardinality $m$ without repetitions, 
then the corresponding quantum minor $D^K_L$ is the element of $\mathcal M$
satisfying
$$
D^K_L = \sum_{\sigma\in\Sigma(m)} \epsilon_P(\sigma K)  
T^{k_{\sigma 1}}_{l_1}\cdots T^{k_{\sigma m}}_{l_m}
= \sum_{\sigma\in\Sigma(m)} \epsilon_P(\sigma K)\epsilon^{-1}_P(\tau L)  
T^{k_{\sigma 1}}_{l_{\tau 1}}\cdots T^{k_{\sigma m}}_{l_{\tau m}}
$$
$$
D^K_L = \sum_{\sigma\in\Sigma(m)} \epsilon_Q(\sigma L)  
T_{l_{\sigma 1}}^{k_1}\cdots T^{l_{\sigma m}}_{k_m}
= \sum_{\sigma\in\Sigma(m)} \epsilon_Q(\sigma L)\epsilon^{-1}_Q(\tau K)  
T_{l_{\sigma 1}}^{k_{\tau 1}}\cdots T_{l_{\sigma m}}^{k_{\tau m}}
$$
where $\tau\in\Sigma(m)$ is a fixed permutation.

\nxpoint {\bf Proposition.}
{\it Let $\iota$ be the monomorphism
~(\ref{eq:iota_SLq_to_multipar}). Then $\iota(D^K_J) = l^J d^K_J r_K$
where $d^K_J$ denotes the quantum minor in 1-parametric case.

Proof.} \[\begin{array}{lcl}
\iota(D^K_J) &=& \sum_{\sigma \in \Sigma(m)}
(-q)^{l(\sigma)} \zeta_r(\sigma K)
\iota(T^{k\sigma(1)}_{j_1})\cdots \iota(T^{k\sigma(m)}_{j_m})
\\&=& \sum_{\sigma \in \Sigma(m)} (-q)^{l(\sigma)}
l^{j_1} l^{j_2} \cdots l^{j_m} \otimes
t^{k\sigma(1)}_{j_1}\cdots  t^{k\sigma(m)}_{j_m}
\otimes \zeta_r(\sigma K) r_{k_1}\cdots r_{k_m}
\\&=& \sum_{\sigma \in \Sigma(m)}
(-q)^{l(\sigma)} l^{j_1} l^{j_2} \cdots l^{j_m} \otimes
t^{k\sigma(1)}_{j_1}\cdots  t^{k\sigma(m)}_{j_m}\otimes r_{k\sigma(1)}
\cdots r_{k\sigma(m)} 
\\&=&  l^J D^K_J r_K.
\end{array}\]

This proposition reminds but is different to the statement of
Lemma 5 in~\cite{Artin:multi} which asserts that the twisting 
considered as an identity map but changing its algebra structure,
interchanges the quantum determinants. Notice that the $D^K_L$ and $d^K_L$
on the two sides are given by different formulas in terms of generators
$T^k_l$ and $t^k_l$ respectively. 

\nxpoint Now one needs to see what happens when one considers
monomials in $D$-s, for example $\iota(D^K_L D^{M}_{N})$.
By the same, method, one gets $(l^L l^M) (d^K_L d^M_N) (r_K r_M)$.
One knows that the relations in $\mathcal M$ are homogeneous
in the sense that the total row multilabel and column multilabel
are the same up to the ordering. Every relation is a sum of
homogeneous. Now if we take different monomials 
$\prod d^*_* = d^{K_1}_{L_1} \cdots d^{K_m}_{L_m}$ in usual quantum minors
then in order to make them manifestly in the image of $\iota$ on
some monomial in multiparametric quantum minors, we need to homogenize
expression by multiplying it by $l^S$ and $r_V$ where $S$ and $V$ are the
ascendingly ordered column and row total multilabel of $\prod d^*_*$, that is
$S = :K_1 \cdots K_m :$ and $V = :L_1\cdots L_m:$. {\it Thus the
$S$ and $V$ are the same for all monomials in the identity (this is more or
less the definition of a homogeneous identity).}
Then  we reorder the multilabels in $l$ and in $r$ separately to get the 
same ordering, but this involves introducing inverse of $\zeta_l$ and $\zeta_r$
corresponding to the ordering on the column and row multilabels
seperately. For a homogeneous identity the multiplier $l^S$ and $r_V$ will
be the same, however the reordering factors will be clearly
different. Thus we get, in terms of $D$-s the same identity up to
different homogeneous factors in front of 
different monomials will be
$\zeta_l^{-1}({L_1}{L_2}\cdots{L_m})\zeta_r^{-1}({K_1}{K_2}\cdots{K_m})$
where $K_i$-s is the row multilabel of $i$-th row and $L_i$ 
of $i$-th column.

{\bf Theorem.} {\it This procedure induces the 1-1 correspondence
between the quantum minor identities for  
the 1-parametric and the minor identities for multiparametric minors.
}

Notice that despite the fact that $\iota$ is not an algebra homomorphism,
essentially $\iota$ and extracting the proprotionality constants from
reorderings of $r$-s and $l$-s do the job. If one would use the
original twisting of \cite{Artin:multi} one has there
a {\it coalgebra map}, hence it sends identities to something
what are not, hence it is not clear how to 
directly use it for the same result.

\end{document}